\documentstyle[12pt,amssymb,amscd]{amsart}  

\newtheorem{thm}{Theorem}[subsection]

\newtheorem{cor}[thm]{Corollary}
\newtheorem{prop}[thm]{Proposition}
\newtheorem{lemma}[thm]{Lemma}

\theoremstyle{remark}
\newtheorem{remark}[thm]{Remark}

\theoremstyle{definition}

\numberwithin{equation}{section}

\newcommand{\bbC}{\Bbb C}
\newcommand{\bbZ}{\Bbb Z}
\newcommand{\bbP}{\Bbb P}
\newcommand{\bbQ}{\Bbb Q}
\newcommand{\bbL}{\Bbb L}

\newcommand{\cO}{\cal O}
\newcommand{\cD}{\cal D}
\newcommand{\cL}{\cal L}
\newcommand{\cM}{\cal M}
\newcommand{\cX}{\cal X}

\newcommand{\Pic}{\operatorname{Pic}}
\newcommand{\divis}{\operatorname{div}}
\newcommand{\Div}{\operatorname{Div}}
\newcommand{\Cl}{\operatorname{Cl}}
\newcommand{\ord}{\operatorname{ord}}
\newcommand{\codim}{\operatorname{codim}}
\newcommand{\D}{\operatorname{D}}
\newcommand{\Perv}{\operatorname{Perv}}
\newcommand{\RH}{\operatorname{RH}}
\newcommand{\DR}{\operatorname{DR}}
\newcommand{\SS}{\operatorname{SS}}
\newcommand{\supp}{\operatorname{Supp}}
\newcommand{\R}{\operatorname{\bold R}}

\newcommand{\pr}{\operatorname{pr}}

\newcommand{\Coker}{\operatorname{Coker}}

\newcommand{\isomoto}{\overset{\sim}{\to}}
\newcommand{\isomo}{\overset{\sim}{=}}

\thanks
{The second author was supported
in part by NSF grant DMS-9504522.}

\begin{document}

\title{On the singularities of Theta divisors on Jacobians}
\author{P.Bressler, J.-L. Brylinski}
\maketitle

\section{Introduction}\label{section:intro}
The theta divisor $\Theta$ of the Jacobian
variety of a complex curve $X$ is best viewed
as a  divisor inside the component $\Pic^{g-1}(X)$
consisting of (isomorphism classes of)
line bundles of degree $g-1$.
Then a line bundle $L$ belongs to $\Theta$
if and only if it has non-zero sections.
Riemann proved that the multiplicity
of $\Theta$ at a point $L$ is equal
to $\dim H^0(X;L)-1$. Kempf (\cite{Ke}) obtained a geometric
proof of Riemann's theorem and a beautiful
description of the tangent cone to $\Theta$
at any point.

In this paper we study the intersection
cohomology of $\Theta$ when $X$ is not
hyperelliptic. Our starting point is a theorem
of Martens concerning the geometry of the
Abel-Jacobi mapping $\phi:S^{g-1}(X)\to \Pic^{g-1}(X)$
and of its fibers. We interpret this theorem
as saying that $\phi$ is small in the sense
of Goresky and MacPherson. This means that the 
intersection cohomology $IH^\bullet(\Theta;\bbQ)$ is isomorphic
to the cohomology $H^\bullet(S^{g-1}(X);\bbQ)$. The cohomology
of $S^{g-1}(X)$, including the algebra structure,
was completely determined by MacDonald in \cite{Mac} . 

From the evaluation of the differential of $\phi$
we deduce (Theorem \ref{thm:main}) that the intersection complex has the 
property that its characteristic variety (inside
the cotangent bundle of $\Pic^{g-1}(X)$) is irreducible.
This is a rather unusual phenomenon; it is known to
be true for Schubert varieties in classical grassmannians
(\cite{BFL}) and more generally in hermitian symmetric
spaces of simply-laced groups (\cite{BF}).

In Section \ref{section:invol} we study the effect on the intersection 
cohomology of the involution $\iota$ of $\Theta$ given by
\[
\iota(L)=\Omega^1_X\otimes L^{\otimes -1}\ .
\]
We show that, for curves of even genus, the action of $\iota$
on $IH^\bullet(\Theta;\bbQ)$ does not preserve the algebra structure
of $H^\bullet(S^{g-1}(X);\bbQ)$. This gives another interpretation
of the well-known ``calcul triste'' of Verdier (\cite{BG}), leading
to more examples of singular varieties with two different
small resolutions yielding different algebra structures
on intersection cohomology.

In Section \ref{section:invol} interpret the classical computation of
the number of odd $\Theta$-characteristics, which is equal to
$2^{2g-1}-2^{g-1}$ according to \cite{W} and \cite{Mu1}, in terms of
intersection cohomology.

A $\Theta$-characteristic $L$ is a square root
of the canonical bundle $\Omega^1_X$. Thus, a $\Theta$-characteristic with
$\dim H^0(X;L) > 0$ determines a fixed point of
$\iota$ in $\Theta$. A $\Theta$-characteristic $L$ is called odd or even 
depending on the parity of $\dim H^0(X;L)$.

Our idea is to apply the Lefschetz fixed point formula to the action of
$\iota$ on $IH^\bullet(\Theta;\bbQ)$. We show that the contribution of
the theta characteristic $L$
to the fixed point formula is equal to $1$ if $L$ is odd
and $0$ if $L$ is even. On the other hand we show that
the (super)trace of $\iota$ acting on intersection
cohomology is equal to $2^{2g-1}-2^{g-1}$.

We would like to thank Alberto Collino for communicating his calculations
to us.

\section{Algebraic curves and their Jacobians}\label{section:jacs}
In what follows $X$ will denote a connected smooth projective algebraic curve
of genus $g$ over the field $\bbC$ of complex numbers, i.e. a compact Riemann
surface.

\subsection{Line bundles}
Let $\Pic(X)$ denote the set of isomorphism classes of algebraic (equivalently
holomorphic) line bundles on $X$. The operation of tensor product of line
bundles endows $\Pic(X)$ with the structure of an algebraic group. The identity 
element is given by the (isomorphism class of) the structure sheaf $\cO_X$.
The map
\begin{eqnarray*}
\deg : \Pic(X) & \to & \bbZ \\
L & \mapsto & \deg(L) = \int_X c_1(L)
\end{eqnarray*}
induces an isomorphism on the respective groups of connected components
\[
\pi_0(\Pic(X))\isomoto\bbZ\ .
\]
As is usual, we denote by $\Pic^d(X)$ the component of $\Pic(X)$ consisting
of (the isomorphism classes of) line bundles of degree $d$.

The connected component of the identity $\Pic^0(X)$ is an Abelian variety
isomorphic to the quotient
$H^1(X;\cO_X)/H^1(X;\bbZ)\isomo 
H^1(X;\bbC)/\left(H^0(X;\Omega^1_X)+H^1(X;\bbZ)\right)$.
In particular the Lie algebra of $\Pic^0(X)$ is $H^1(X;\cO_X)$, and the
cotangent space at the identity is $H^0(X;\Omega^1_X)$. The component
$\Pic^d(X)$ is a principal homogeneous space under $\Pic^0(X)$.

There is an involution
\begin{eqnarray*}
\iota : \Pic(X) & \to & \Pic(X) \\
L & \mapsto & \Omega^1_X\otimes L^{\otimes -1}
\end{eqnarray*}
which maps $\Pic^d(X)$ to $\Pic^{2g-2-d}(X)$ and, therefore preserves
$\Pic^{g-1}(X)$.

\subsection{Divisors}
Let $\Div(X)$ denote the free Abelian group generated by the (closed) points
of $X$. An element $D = \sum_i m_i\cdot p_i$ (where $m_i\in\bbZ$ and
$p_i\in X$, and the sum is finite) is called a {\em divisor} on $X$.
The divisor $D$ is said to be {\em effective} if $m_i\geq 0$ for all $i$.
We will denote this fact by $D\geq 0$. The {\em degree} of the divisor $D$ is
the integer defined by $\deg(D) = \sum_i m_i$.

To a nonzero meromorphic function $f$ on $X$ one associates the divisor
$\divis(f) = \sum_{p\in X} \ord_p(f)\cdot p$. Note that $\deg(\divis(f)) = 0$.
The divisors of the form
$\divis(f)$ are called {\em principal} and form the subgroup $P(X)$ of
$\Div(X)$. The quotient group $\Div(X)/P(X)$ is denoted by $\Cl(X)$ and
is called the {\em divisor class group} of $X$.

To a nonzero meromorphic section $s$ of a line bundle $L$ on $X$ one
associates the divisor $\divis(s) = \sum_{p\in X} \ord_p(s)\cdot p$.
If, in addition, $f$ is a nonzero meromorphic function on $X$, then
$\divis(fs) = \divis(f) + \divis(s)$. Since for any two nonzero meromorphic
sections $s_1$ and $s_2$ of $L$ one has a (necessarily nonzero) meromorphic
function $f$ such that $s_1 = fs_2$, the divisors $\divis(s_1)$ and $\divis(s_2)$ are in the same divisor class in $\Cl(X)$.

Moreover, if $s_1$ and $s_2$ are nonzero meromorphic sections of line bundles
$L_1$ and $L_2$ respectively, then $s_1\otimes s_2$ is a nonzero meromorphic
section of $L_1\otimes L_2$ and $\divis(s_1\otimes s_2) = \divis(s_1) + \divis(s_2)$. 

Thus, the association $L\mapsto \divis(s)$, where $s$ is a nonzero meromorphic 
section of $L$, gives rise to a well defined homomorphism
\begin{eqnarray*}
\Pic(X) & \to & \Cl(X) \\
L & \mapsto & D(L)\ .
\end{eqnarray*}
This map is in fact an isomorphism. The inverse is given by the map
\begin{eqnarray*}
\Cl(X) & \to & \Pic(X) \\
D & \mapsto & \cO_X(D)\ ,
\end{eqnarray*}
where $\cO_X(D)$ is the subsheaf of the sheaf of meromorphic functions on $X$
consisting of functions $f$ which are holomorphic in the complement of $D$ and
satisfy $\divis(f)\geq D$. In particular, if $D$ is effective the line bundle
$\cO_X(D)$ has nonzero holomorphic sections, i.e. $H^0(X;L)\neq 0$.

\subsection{Effective divisors}
The set of effective divisors of degree $d$ on $X$ is easily identified with
the $d$-th symmetric power $S^d(X)$ of the curve $X$. The variety
$S^d(X) = X^{\times d}/\Sigma_d$ is smooth of dimension $d$ for all $d\geq 0$.
For an effective divisor $D$, the {\em complete linear system} $\vert D\vert$
is defined as the set of all effective divisors in the class of $D$ in $\Cl(X)$.
The association $s\mapsto \divis(s)$ for a nonzero section
$s\in H^0(X;\cO_X(D))$ gives rise to the natural isomorphism 
$\bbP(H^0(X;L))\isomoto\vert D\vert$.

The map
\begin{eqnarray*}
\phi : S^d(X) & \to & \Pic^d(X) \\
D & \mapsto & \cO_X(D)
\end{eqnarray*}
is a morphism of algebraic varieties. It is surjective for $d\geq g$ and
birational onto its image for $d\leq g$. For $L\in\Pic^d(X)$ the fiber
$\phi^{-1}(L)$ is naturally identified with the complete linear system
$\bbP(H^0(X;L))$.

\subsection{The $\Theta$-divisor}
The $\Theta$-divisor is defined as the image of the map
\[
\phi : S^{g-1}(X) \to \Pic^{g-1}(X)
\]
and will be denoted by $\Theta$. It is an irreducible closed subvariety
of $\Pic^{g-1}(X)$ (for $X$ smooth connected) of codimension one. From
the discussion above it follows that $\Theta$ is the locus of (isomorphism
classes of) line bundles $L$ of degree $g-1$ with $H^0(X;L)\neq 0$.

The hypersurface $\Theta$ is singular in general. It is known that
$\dim Sing(\Theta)\geq g-4$. The singular locus is
determined from the following result of Riemann (see \cite{Ke}).

\begin{thm}
The multiplicity of $\Theta$ at the point $L$ is equal to
$\dim H^0(X;L)-1$.
\end{thm}
\begin{cor}
A line bundle $L$ of degree $g-1$ determines a singular point of $\Theta$
if and only if $\dim H^0(X;L)\geq 2$.
\end{cor}

For a line bundle $L$ of degree $g-1$ the Riemann-Roch theorem shows that
\linebreak
$\dim H^0(X;L) = \dim H^0(X;\Omega^1_X\otimes L^{\otimes -1})$. Therefore
the $\Theta$-divisor is preserved by the involution $\iota$.

The $\Theta$-divisor is naturally stratified by the closed subvarieties
$W^r_{g-1}$ defined as the locus of lined bundles $L$ of degree $g-1$
with $\dim H^0(X;L) - 1\geq r$. In particular $\Theta = W^0_{g-1}$ and
$W^1_{g-1}$ is the singular locus of $\Theta$. The following theorem of
Martens, stated here in the particular case of interest, gives an estimate on 
the dimension of $W^r_{g-1}$. Note that, by Clifford's theorem,
$2\cdot\dim H^0(X;L)\leq g-1$.

\begin{thm}
Suppose that $g\geq 3$, $X$ not hyperelliptic and $2r\leq g-1$.
Then all components of $W^r_{g-1}$ have the same dimension, and
$\dim W^r_{g-1}\leq g - 2r - 2$.
\end{thm}

Recall that a map $f: Y\to Z$ of algebraic varieties is called {\em small}
if
\[
\codim \lbrace z\in Z\vert\ \dim f^{-1}(z)\geq d\rbrace > 2d\ .
\]
Thus Martens' theorem has the following corollary.

\begin{cor}\label{cor:small}
Suppose that $X$ is not hyperelliptic. Then the map
\linebreak $\phi : S^{g-1}(X)\to\Theta$ is small.
\end{cor}

\section{The characteristic variety of the intersection complex}
\label{section:char-var}
From now on we will assume that $X$ is not hyperelliptic.

\subsection{The Riemann-Hilber correspondence}
For an algebraic variety $Y$ over $\bbC$ let $\D^b_c(Y;k)$ denote the
bounded derived category of complexes of sheaves of $k$-vector spaces on
$Y(\bbC)$ with (algebraically) constructible cohomology. Let $\Perv(Y;k)$
denote the full (Abelian) subcategory of perverse sheaves.

For a smooth algebraic variety $Y$ over $\bbC$ let $\D^b_{rh}(\cD_Y)$
denote the bounded derived category of complexes of left $\cD_Y$-modules
with regular holonomic cohomology. Let $\RH(\cD_Y)$ denote the full
(Abelian) subcategory of left holonomic $\cD_Y$-modules with regular 
singularities.

Recall that the de Rham functor
\begin{eqnarray*}
\DR : \D^b_{rh}(\cD_Y) & \to & \D^b_c(Y;\bbC) \\
M^\bullet & \mapsto & \omega_Y\otimes^{\bbL}_{\cD_Y}M^\bullet[\dim Y]
\end{eqnarray*}
is an equivalence of categories called the {\em Riemann-Hilbert
correspondence} (\cite{KK}, \cite{Me1}, \cite{Me2}, see also \cite{B})),
and restricts to the exact equivalence of Abelian categories
\[
\DR : \RH(\cD_Y)\to\Perv(Y;\bbC)\ .
\]

\subsection{The intersection cohomology of the $\Theta$-divisor}
Let $\Theta^{reg} = \Theta\setminus W^1_{g-1}$ denote the nonsingular
part of the $\Theta$-divisor and let
$i : \Theta^{reg}\hookrightarrow \Pic^{g-1}(X)$ denote the (locally closed)
inclusion map. From Corollary \ref{cor:small} and \cite{GM1} we obtain the
following proposition.

\begin{prop}\label{prop:IC-is-dir-im}
There is an isomorphism
$\R\phi_*(\bbQ_{S^{g-1}(X)})\isomo i_{!*}(\bbQ_{\Theta^{reg}})$ in 
\linebreak $\D^b_c(\Pic^{g-1}(X);\bbQ)$  
\end{prop}

Here $i_{!*}(\bbQ_{\Theta^{reg}})$ is the ``middle'' (\cite{GM1}, \cite{BBD}) 
extension
of the constant sheaf $\bbQ_{\Theta^{reg}}$ such that, in particular,
there is an isomorphism $H^\bullet(\Pic^{g-1}(X);i_{!*}(\bbQ_{\Theta^{reg}}))
\isomo IH^\bullet(\Theta;\bbQ)$. It follows from Proposition 
\ref{prop:IC-is-dir-im} that there is an isomorphism
\[
IH^\bullet(\Theta;\bbQ)\isomo H^\bullet(S^{g-1}(X);\bbQ)\ .
\]

The $\cD_{S^{g-1}(X)}$-module associated (under the Riemann-Hilbert 
correspondence) to the constant sheaf $\bbC_{S^{g-1}(X)}$ is the sructure
sheaf $\cO_{S^{g-1}(X)}$. The Corollary \ref{cor:small} implies that
the cohomology of the direct image (in $\cD$-modules) $\phi_+\cO_{S^{g-1}(X)}$
is a complex with cohomology concentrated only in degree zero so that
$\phi_+\cO_{S^{g-1}(X)}\isomo H^0\phi_+\cO_{S^{g-1}(X)}$ in
$\D^b_{rh}(\cD_{\Pic^{g-1}(X)})$. Let
$\cL$ denote the regular holonomic $\cD_{\Pic^{g-1}(X)}$-module
$H^0\phi_+\cO_{S^{g-1}(X)}$, so that $\DR(\cL)\isomo i_{!*}\bbC_{\Theta^{reg}}$.
The $\cD$-module $\cL$ may be characterized
(up to a unique isomorphism) as the smallest nontrivial submodule of the 
$\cD$-module $i_+\cO_{\Theta^{reg}}$ (see \cite{Ka} and also \cite{B}).

\subsection{The characteristic variety}
Recall that to a holonomic $\cD_Y$-module $M$ on a complex algebraic variety
$Y$ one associates the characteristic cycle $\SS(M)$ which is an effective
conic Lagrangian cycle on the cotangent bundle $T^*Y$ by a theorem of
Sato-Kashiwara-Kawai, Malgrange and Gabber (\cite{SKK}, \cite{Mal}, \cite{G}). 
It is known 
(\cite{Ka}, see also \cite{B}) that it is of the form
\[
\SS(M) = \sum_i m_i\cdot \overline{T^*_{Y_i}Y}
\]
for suitable smooth locally closed subvarieties $Y_i$ of $\supp M$ and
positive integers $m_i$. For example, if $f:Z\hookrightarrow Y$
is the
inclusion of a closed smooth subvariety $Z$, then $\SS(f_+\cO_Z) = T^*_ZY$.
The multiplicities of the components of the characteristic cycle are local
in the sence that $m_i$ depends only on the restriction of $M$ to any open
(in $Y$) neighborhood of any point of $Y_i$.

From the discussion above we may conclude that the characteristic cycle of
$\cL$ is of the form
\[
\SS(\cL) = \overline{T^*_{\Theta^{reg}}\Pic^{g-1}(X)} +
\sum_i m_i\cdot\overline{T^*_{Y_i}\Pic^{g-1}(X)}
\]
for suitable smooth locally closed subvarieties $Y_i$ of $\Theta$.

The main result of this note is the following.

\begin{thm}\label{thm:main}
Suppose that $X$ is a smooth connected projective algebraic curve over
$\bbC$. Let $\cL$ denote the simple holonomic $\cD_{\Pic^{g-1}(X)}$-module
with regular singularities which restricts to $\cO_{\Theta^{reg}}$ on the
nonsingular part $\Theta^{reg}$ of the $\Theta$-divisor.

Then the characteristic cycle of $\cL$ is irreducible,
i.e. $\SS(\cL) = \overline{T^*_{\Theta^{reg}}\Pic^{g-1}(X)}$. 
\end{thm}

\begin{pf}
We have the isomorphism (in the derived category) $\phi_+\cO_{S^{g-1}(X)}\isomo
\cL$ and $\SS(\cO_{S^{g-1}(X)}) = T^*_{S^{g-1}(X)}{S^{g-1}(X)}$. Thus,
according to Kashiwara, there is an inclusion
\[
\supp\SS(\cL)\subseteq\pr((d\phi^t)^{-1}(T^*_{S^{g-1}(X)}{S^{g-1}(X)}))
\]
where the maps
\[
T^*S^{g-1}(X) @<{d\phi^t}<< S^{g-1}(X)\times_{\Pic^{g-1}(X)}
T^*\Pic^{g-1}(X) @>{\pr}>> T^*\Pic^{g-1}(X)
\]
are the (transpose of) the differential of the map $\phi$ and the projection
on the second factor. We have the following classical identification.

\begin{lemma}
$\ker(d\phi^t)\isomo H^0(X;\Omega^1_X\otimes\cO_X(-D))\subset
H^0(X;\Omega^1_X)$
\end{lemma}

Therefore the subvariety 
\[
(d\phi^t)^{-1}(T^*_{S^{g-1}(X)}{S^{g-1}(X)}) =
\lbrace (D,\xi)\vert\ D\in S^{g-1}(X),\ \xi\in T^*_{\phi(D)}\Pic^{g-1}(X),
\ d\phi^t(\xi) = 0\rbrace
\]
is naturally described as
\[
(d\phi^t)^{-1}(T^*_{S^{g-1}(X)}{S^{g-1}(X)}) =
\lbrace (D,\omega)\vert\ D\in S^{g-1}(X),\ \omega\in 
H^0(X;\Omega^1_X\otimes\cO_X(-D))\rbrace\ .
\]
From this description one sees immediately that 
$(d\phi^t)^{-1}(T^*_{S^{g-1}(X)}{S^{g-1}(X)})$ is the union of irreducible
components $\overline Z^r$, where the locally closed subvariety $Z^r$ of
\linebreak $S^{g-1}(X)\times_{\Pic^{g-1}(X)}T^*\Pic^{g-1}(X)$ is given by
\[
Z^r = \lbrace (D,\omega)\vert\ \phi(D)\in W^r_{g-1}\setminus W^{r+1}_{g-1},
\ \omega\in H^0(X;\Omega^1_X\otimes\cO_X(-D))\rbrace\ .
\]
In particular $\dim Z^r = \dim W^r_{g-1} + \dim\bbP(H^0(X;L)) +
\dim H^0(X;\Omega^1_X\otimes L^{\otimes -1})$, where $L$ is a general point
of $W^r_{g-1}$. From the theorems of Clifford and Martens one deduces 
immediately that $\dim Z^0 = g$ and $\dim Z^r\leq g-1$ for $r>0$.
Therefore $\dim\pr(\overline Z^r)\leq g-1$ for $r>0$. It is easy to see that
$\pr(Z^0) = \overline{T^*_{\Theta^{reg}}\Pic^{g-1}(X)}$. 

Since all components of $\supp\SS(\cL)$ are Lagrangian and, therefore, of 
dimension exactly $g$, one must have the inclusion 
$\supp\SS(\cL)\subseteq\pr(Z^0)$.
\end{pf}

We also obtain the following description of
$\overline{T^*_{\Theta^{reg}}\Pic^{g-1}(X)}$: it is a conical Lagrangian
subvariety of $T^*\Pic^{g-1}(X)$ which projects to $\Theta$; the fiber
over $L$ is given by
\[
\overline{T^*_{\Theta^{reg}}\Pic^{g-1}(X)}\cap T^*_L\Pic^{g-1}(X)
= \bigcup_{D\in\bbP(H^0(X;L))}H^0(X;\Omega^1_X(-D))\ ,
\]
where we use the natural identification $T^*_L\Pic^{g-1}(X)\isomo
H^0(X;\Omega^1_X)$.

\begin{remark}
It follows from Theorem \ref{thm:main} and a theorem of Andronikof (\cite{A})
that the wave-front set of the current defined by integrating over $\Theta$
is irreducible.
\end{remark}

\subsection{The universal $\Theta$-divisor}
Let $\cM_g^{(n)}$ denote the moduli space of curves of genus $g$ with
level $n$ structure. Let $\cX\to\cM_g^{(n)}$ denote the universal curve of
genus $g$. One has the varieties $\Pic^d(\cX/\cM_g^{(n)})$ for $d\in\bbZ$,
the universal $\Theta$-divisor 
$\Theta_{univ}\hookrightarrow\Pic^{g-1}(\cX/\cM_g^{(n)})$ and
the map $\phi : S^{g-1}(\cX/\cM_g^{(n)})\to\Pic^{g-1}(\cX/\cM_g^{(n)})$
of varieties over $\cM_g^{(n)}$ which is birational onto $\Theta_{univ}$.
In this setting we have the analog of Theorem \ref{thm:main}

\begin{thm}
The map $\phi : S^{g-1}(\cX/\cM_g^{(n)})\to\Pic^{g-1}(\cX/\cM_g^{(n)})$
is small. Consequently there is a canonical isomorphism (of
$H^\bullet(\cM_g^{(n)};\bbQ)$-modules)\linebreak
$IH^\bullet(\Theta_{univ};\bbQ)\isomo
H^\bullet(S^{g-1}(\cX/\cM_g^{(n)});\bbQ)$.

The characteristic variety (cycle) of the simple
$\cD_{\Pic^{g-1}(\cX/\cM_g^{(n)})}$-module
which restricts to $\cO_{\Theta_{univ}^{reg}}$ on the nonsingular part
$\Theta_{univ}^{reg}$ of $\Theta_{univ}$ is irreducible.
\end{thm}

\section{Action of the involution}\label{section:invol}
In this section we discuss the induced action of the involution
\begin{eqnarray*}
\iota: \Pic^{g-1}(X) & @>>> & \Pic^{g-1}(X) \\
L & \mapsto & \Omega^1_X\otimes L^{\otimes -1}
\end{eqnarray*}
on the intersection cohomology of the $\Theta$-divisor.

By Corollary \ref{cor:small} there is an isomorphism
\[
H^\bullet(S^{g-1}(X);\bbQ)\isomo IH^\bullet(\Theta;\bbQ)
\]
which induces an algebra structure on $IH^\bullet(\Theta;\bbQ)$.
We will show that the involution
$\iota$ does not preserve this algebra structure.

In what follows we will not make notational distictions between
homology classes of cycles and their Poincare duals in cohomology.

\subsection{The Riemann-Roch correspondence}
Let $\rho\in H^{2g-2}(S^{g-1}(X)\times S^{g-1}(X);\bbQ)$ denote the
class of the cycle
\[
\left\lbrace (D_1,D_2)\in S^{g-1}(X)\times S^{g-1}(X)\vert
D_1+D_2\in\vert K\vert\right\rbrace\ ,
\]
where $K$ denotes the canonical divisor, which we call
{\em the Riemann-Roch correspondence}.

The canonical linear system
\[
\vert K\vert\isomo \bbP(H^0(X;\Omega^1_X))\hookrightarrow S^{2g-2}(X)
\]
determines the class $\kappa\in H^{2g-2}(S^{2g-2}(X);\bbQ)$. The class
$\rho$ is the image of the class $\kappa$ under the pullback map
\[
\Sigma^* : H^{2g-2}(S^{2g-2}(X);\bbQ) @>>> H^{2g-2}(S^{g-1}(X)\times
S^{g-1};\bbQ)
\]
induced by the map
\begin{eqnarray*}
\Sigma : S^{g-1}(X)\times S^{g-1}(X) & @>>> & S^{2g-2}(X) \\
(D_1,D_2) & \mapsto & D_1+D_2 \ .
\end{eqnarray*}

The Riemann-Roch correspondence acts on $H^\bullet(S^{g-1}(X);\bbQ)$
by
\[
\alpha\mapsto (\pr_2)_*(pr_1^*(\alpha)\smile\rho)\ ,
\]
where $\pr_i: S^{g-1}(X)\times S^{g-1}(X)\to S^{g-1}(X)$ denotes
the projection on the $i^{\text{th}}$ factor.

\begin{prop}
Under the isomorphism $IH^\bullet(\Theta;\bbQ)\isomo H^\bullet(S^{g-1}(X);
\bbQ)$ the action of the involution $\iota$ on $IH^\bullet(\Theta;\bbQ)$
corresponds to the action of the Riemann-Roch correspondence on
$H^\bullet(S^{g-1}(X);\bbQ)$.
\end{prop}

\begin{remark}
It follows, in particular, that the action of the Riemann-Roch correspondence
is, in fact, an involution of $H^\bullet(S^{g-1}(X);\bbQ)$.\qed
\end{remark}

\subsection{Cohomology of symmetric powers of a curve}
A complete study of the cohomology of a symmetric power of a curve may
be found in \cite{Mac}. We will need the following facts.

The Abel-Jacobi map
\[
\phi_d : S^d(X) @>>> \Pic^d(X)
\]
induces on $H^\bullet(S^d(X);\bbQ)$ the structure of a module over
$H^\bullet(\Pic^d(X);\bbQ)$. 

A point $p\in X$ determines the embedding
\begin{eqnarray*}
j : S^{d-1} & @>>> & S^d(X) \\
D & \mapsto & D+p
\end{eqnarray*}
and the isomorphism
\[
\otimes\cO(p) : \Pic^{d-1}(X) @>>> Pic^d(X)
\]
such that $\phi_d\circ j= (\otimes\cO(p))\circ\phi_{d-1}$. The homotopy
classes of the maps $j$ and $\otimes\cO(p)$ do not depend on the point $p$
for $X$ connected. In particular the induced map 
$j^* : H^\bullet(S^d(X);\bbQ)\to H^\bullet(S^{d-1}(X);\bbQ)$ and the
Gysin map $j_* : H^\bullet(S^{d-1}(X);\bbQ)\to H^{\bullet +2}(S^d(X);\bbQ)$
are well defined.

The map $(\otimes\cO(p))^*$ provides the canonical
identification of algebras $H^\bullet(\Pic^d(X);\bbQ)$ for various $d$. 
Note that there is a canonical isomorphism
$H^\bullet(\Pic^d(X);\bbQ)\isomo{\bigwedge}^\bullet H^1(X;\bbQ)$.

Moreover, $j^*$ and $j_*$ are maps of modules over
${\bigwedge}^\bullet H^1(X;\bbQ)$, $j^*$ is surjective and $j_*$ is
injective.

Let $\eta_d = j_*(1) = j_*([S^{d-1}(X)])\in H^2(S^d(X);\bbQ)$. Then
the identity $j^*(\eta_d) = \eta_{d-1}$ holds, hence $j_*(\eta_{d-1})=
\eta_d^2$

Consider the map
\begin{eqnarray*}
\Sigma : S^{d_1}(X)\times S^{d_2}(X) & @>>> & S^{d_1+d_2}(X) \\
(D_1,D_2) & \mapsto & D_1+D_2\ .
\end{eqnarray*}
Then $\Sigma^*(\eta_{d_1+d_2}) = \pr_1^*(\eta_{d_1})+\pr_2^*(\eta_{d_2})$.

Multiplication by $\eta_d$ commutes with the action of
${\bigwedge}^\bullet H^1(X;\bbQ)$. Moreover, the natural map
\[
{\bigwedge}^\bullet H^1(X;\bbQ)\otimes\bbQ [\eta] @>>> H^\bullet(S^d(X);
\bbQ)
\]
is surjective in all degrees and an isomorphism in degrees up through $d$.

The class $\theta = [\Theta]\in H^2(\Pic^{g-1}(X);\bbQ)\isomo{\bigwedge}^2
H^1(X;\bbQ)$ corresponds to the (symplectic) intersection pairing on $X$.
The Poincar\'e's formula (\cite{ACGH}) says that
\[
(\phi_{g-m})_*(1) = (\phi_{g-m+i})_*(\eta_{g-m+i}^i) =
\frac{\theta^m}{m!}\ .
\]
In particular, for $m=g$ we find that $\displaystyle\frac{\theta^g}{g!}$
is the class of a point. 

The map $\phi_{2g-1}:S^{2g-1}(X)\to\Pic^{2g-1}(X)$ is a projective space
bundle.
The map $j:S^{2g-2}(X)\hookrightarrow S^{2g-1}(X)$ restricts to an embedding
of $\vert K\vert$ as a fiber of $\phi_{2g-1}$. Therefore
the identities $j_*(\kappa) = \phi_{2g-1}^*(\displaystyle\frac{\theta^g}{g!})$,
and, consequently, $\phi_{2g-2}^*(\displaystyle\frac{\theta^g}{g!}) =
j^*\phi_{2g-1}^*(\displaystyle\frac{\theta^g}{g!}) =\eta_{2g-2}\smile\kappa$
hold.

\subsection{The action of the Riemann-Roch correspondence}
The action of the involution $\iota$ on $H^p(\Pic^{g-1}(X);\bbQ)$
is given by multiplication by $(-1)^p$ and may be realized as the
action of a correspondence as follows.

Let $\Delta^-\in H^{2g}(\Pic^{g-1}(X)\times\Pic^{g-1}(X);\bbQ)$ denote
the ``antidiagonal'', i.e. the cycle
\[
\left\lbrace (L_1,L_2)\in\Pic^{g-1}(X)\times\Pic^{g-1}(X)\vert
L_1\otimes L_2\isomo\Omega^1_X\right\rbrace\ .
\]
The class $\Delta^-$ is the image of the class of the point
$[\Omega^1_X]\in H^{2g}(\Pic^{2g-2}(X);\bbQ)$ under the map
\[
\otimes^*: H^\bullet(\Pic^{2g-2}(X);\bbQ) @>>>
H^\bullet(\Pic^{g-1}\times\Pic^{g-1}(X);\bbQ)
\]
induced by the map
\begin{eqnarray*}
\otimes : \Pic^{g-1}(X)\times\Pic^{g-1}(X) & @>>> & \Pic^{2g-2}(X) \\
(L_1, L_2) & \mapsto & L_1\otimes L_2\ .
\end{eqnarray*}

Then, clearly, the action of the involution $\iota$ is given by
\[
\iota(\alpha) = (\pr_2^\prime)_*((\pr_1^\prime)^*(\alpha)\smile\Delta^-)\ ,
\]
where $\pr_i^\prime : \Pic^{g-1}(X)\times\Pic^{g-1}(X)\to\Pic^{g-1}(X)$
denotes the projection on the $i^{\text{th}}$ factor.

The action of the Riemann-Roch correspondence on $H^\bullet(S^{g-1}(X);\bbQ)$
is compatible with the action of $\iota$ on $H^\bullet(\Pic^{g-1}(X);\bbQ)$
and the action of $H^\bullet(\Pic^{g-1}(X);\bbQ)$ on
$H^\bullet(S^{g-1}(X);\bbQ)$. Thus, the action of the Riemann-Roch
correspondence is determined by its values on the powers of the class
$\eta_{g-1}$.

The following proposition is due to Alberto Collino.

\begin{prop}
\begin{equation}\label{formula:RR-eta}
(\pr_2)_*(\pr_1^*(\eta_{g-1}^m)\smile\rho) =
\sum_{i=0}^{m}\frac{(-1)^i}{(m-i)!}\phi_{g-1}^*(\theta^{m-i})\smile
\eta_{g-1}^i
\end{equation}
\end{prop}
\begin{pf}
Consider the commutative diagram
\[
\begin{CD}
S^{g-1}(X)\times S^{g-1}(X) @>{\Sigma}>> S^{2g-2}(X) \\
@V{\phi_{g-1}\times\phi_{g-1}}VV @VV{\phi_{2g-2}}V \\
\Pic^{g-1}(X)\otimes\Pic^{g-1}(X) @>{\otimes}>> \Pic^{2g-2}(X)\ .
\end{CD}
\]
We have
\begin{eqnarray*}
(\phi_{g-1}\times\phi_{g-1})^*\Delta^- & = & (\phi_{g-1}\times\phi_{g-1})^*
\otimes^*(\frac{\theta^g}{g!}) \\
& = & \Sigma^*\phi_{2g-2}^*(\frac{\theta^g}{g!}) \\
& = & \Sigma^*(\eta_{2g-2}\smile\kappa) \\
& = & (\pr_1^*(\eta_{g-1})+\pr_2^*(\eta_{g-1}))\smile\Sigma^*\kappa \\
& = & \pr_1^*(\eta_{g-1})\smile\rho + \pr_2^*(\eta_{g-1})\smile\rho\ .
\end{eqnarray*}
Therefore, for $\alpha\in H^\bullet(S^{g-1}(X);\bbQ)$, we have
\begin{multline*}
(\pr_2)_*(\pr_1^*(\alpha)\smile(\phi_{g-1}\times\phi_{g-1})^*\Delta^-) \\ =
(\pr_2)_*(\pr_1^*(\alpha\smile\eta_{g-1})\smile\rho) +
(-1)^{\deg\alpha}\cdot\eta_{g-1}\smile(\pr_2)_*(\pr_1^*(\alpha)\smile\rho)\ .
\end{multline*}

On the other hand,
by Lemma \ref{lemma:ind-cor} below (with $X=Y=S^{g-1}(X)$, $A=B=\Pic^{g-1}(X)$,
$f=g=\phi_{g-1}$ and $\Gamma=\Delta^-$), we have
\[
(\pr_2)_*(\pr_1^*(\alpha)\smile(\phi_{g-1}\times\phi_{g-1})^*\Delta^-) = 
\phi_{g-1}^*\iota^*(\phi_{g-1})_*(\alpha) =
(-1)^{\deg\alpha}\phi_{g-1}^*(\phi_{g-1})_*(\alpha)\ .
\]

Thus, putting two calculations together we obtain the identity
\[
(-1)^{\deg\alpha}\phi_{g-1}^*(\phi_{g-1})_*(\alpha) =
(\pr_2)_*(\pr_1^*(\alpha\smile\eta_{g-1})
\smile\rho)+(-1)^{\deg\alpha}\cdot\eta_{g-1}\smile(\pr_2)_*(\pr_1^*(\alpha)
\smile\rho)\ .
\]
In particular, for $\alpha=\eta_{g-1}^{m-1}$ we have
\[
\phi_{g-1}^*\left(\frac{\theta^m}{m!}\right) =
(\pr_2)_*(\pr_1^*(\eta_{g-1}^m)\smile\rho)+
\eta_{g-1}\smile(\pr_2)_*(\pr_1^*(\eta_{g-1}^{m-1})\smile\rho)
\]
or, equivalently,
\begin{equation}\label{formula:rec-rel}
(\pr_2)_*(\pr_1^*(\eta_{g-1}^m)\smile\rho) =
-\eta_{g-1}\smile(\pr_2)_*(\pr_1^*(\eta_{g-1}^{m-1})\smile\rho)
+\phi_{g-1}^*\left(\frac{\theta^m}{m!}\right)\ .
\end{equation}
In particular, for $m=1$ we have
\[
(\pr_2)_*(\pr_1^*(\eta_{g-1})\smile\rho) =
-\eta_{g-1}
+\phi^*(\theta)
\]
which is \eqref{formula:RR-eta} in this case.

Proceding by induction on $m$ with the help of \eqref{formula:rec-rel} we
obtain \eqref{formula:RR-eta}.
\end{pf}

\begin{lemma}\label{lemma:ind-cor}
Suppose given oriented manifolds $X,\ Y,\ A,\ B$,
maps $f:X\to A$ and $g:Y\to B$, and a class $\Gamma\in H^\bullet(A\times B;
\bbQ)$. Then, for any class $\alpha\in H^\bullet(X;\bbQ)$,
\[
(\pr_Y)_*(\pr_X^*(\alpha)\smile (f\times g)^*\Gamma) =
g^*((\pr_B)_*(\pr_A^*(f_*\alpha)\smile\Gamma))\ .
\]
\end{lemma}

Let
\[
STr(\iota, IH^\bullet(\Theta;\bbQ)) =
\sum_p (-1)^p Tr(\iota : IH^p(\Theta;\bbQ)\to IH^p(\Theta;\bbQ))\ .
\]

\begin{prop}\label{prop:STr}
\[
STr(\iota, IH^\bullet(\Theta;\bbQ)) = 2^{2g-1}-2^{g-1}\ .
\]
\end{prop}
\begin{pf}
Using Poincare duality we can write
\begin{multline}
STr(\iota, IH^\bullet(\Theta;\bbQ)) =
\sum_{p < g-1} (-1)^p\cdot 2\cdot Tr(\iota : IH^p(\Theta;\bbQ)\to 
IH^p(\Theta;\bbQ)) \\
+ (-1)^{g-1} Tr(\iota : IH^{g-1}(\Theta;\bbQ)\to IH^{g-1}(\Theta;\bbQ))\ .
\end{multline}

For $p\leq g-1$ we have
\[
IH^p(\Theta;\bbQ) = \bigoplus_{j} H^{p-2j}(\Pic^{g-1}(X);\bbQ)
\smile\eta_{g-1}^j\ .
\]
Formula \eqref{formula:RR-eta} shows that $\iota$ preserves 
the filtration $F_\bullet IH^p(\Theta;\bbQ)$ defined by
\[
F_qIH^p(\Theta;\bbQ) = \bigoplus_{j\leq q}
H^{p-2j}(\Pic^{g-1}(X);\bbQ)\smile\eta_{g-1}^j
\]
and acts by $(-1)^{p-q}$ on $Gr^F_q IH^p(\Theta;\bbQ)$.
Therefore
\begin{multline*}
Tr(\iota : IH^p(\Theta;\bbQ)\to IH^p(\Theta;\bbQ)) = \\
= Tr(Gr^F_\bullet\iota : Gr^F_\bullet IH^p(\Theta;\bbQ)\to Gr^F_\bullet
IH^p(\Theta;\bbQ)) = \\
= \sum_{i+2j = p}(-1)^{i+j}{2g\choose i}
\end{multline*}
and
\begin{multline}
STr(\iota, IH^\bullet(\Theta;\bbQ)) = \\
= \sum_{p < g-1} (-1)^p\cdot 2\cdot
\sum_{i+2j = p}(-1)^{i+j}{2g\choose i}
+ (-1)^{g-1}
\sum_{i+2j = g-1}(-1)^{i+j}{2g\choose i} = \\
= \sum_{i = 0}^{g-1}\delta(g-1-i){2g\choose i}\ ,
\end{multline}
where $\delta:\bbZ/4\cdot\bbZ\to\bbZ$ is defined by
\[
\delta(0) = 1,\ \delta(1) = 2,\ \delta(2) = 1\ \delta(3) = 0\ .
\]
The proposition now follows from Lemma \ref{lemma:calc} below.
\end{pf}

\begin{lemma}\label{lemma:calc}
$\sum_{i = 0}^{g-1}\delta(g-1-i){2g\choose i} = 2^{2g-1}-2^{g-1}$
\end{lemma}
\begin{pf}
First observe that
\[
2^{2g-1} = \frac12 (1+1)^{2g}=
\sum_{i=0}^{g-1}{2g\choose i} +\frac12{2g\choose g}\ .
\]
Now consider separately the four cases corresponding to the residue
of $g$ modulo $4$.

Suppose that $g$ is divisible by $4$. Then
\[
2^{g-1} = \frac12 (1+\sqrt{-1})^{2g} = \sum_{i=0}^{g-1}\frac12
(\sqrt{-1}^i + \sqrt{-1}^{-i}){2g\choose i} +\frac12{2g\choose g}\ .
\]
Therefore
\[
2^{2g-1} - 2^{g-1} = \sum_{i=0}^{g-1}\left( 1 - \frac12
(\sqrt{-1}^i + \sqrt{-1}^{-i})\right){2g\choose i}\ .
\]
Now observe that
\[
 1 - \frac12 (\sqrt{-1}^i + \sqrt{-1}^{-i}) = \delta(-1-i) =
\delta(g-1-i)
\]
and the equality follows.

Other cases follow from similar calculations and we omit the details.
If $g\equiv 2\bmod 4$, then we find
that
$$2^{2g-1}-2^{g-1}={1\over 2}(1+1)^{2g}+{1\over 2}
(1+\sqrt{-1})^{2g}=\sum_{j=1}^{g-1}\delta(1-j){2g\choose j}.$$
If $g\equiv 1\bmod 4$ , then we find
that
$$2^{2g-1}-2^{g-1}={1\over 2}(1+1)^{2g}
+{\sqrt{-1}\over 4}(1+\sqrt{-1})^{2g}
-{\sqrt{-1}\over 4}(1-\sqrt{-1})^{2g}
=\sum_{j=0}^{g-1}\delta(-j){2g\choose j}.$$
If $g\equiv 3\bmod 4$ , then we find
that
$$2^{2g-1}-2^{g-1}={1\over 2}(1+1)^{2g}
-{\sqrt{-1}\over 4}(1+\sqrt{-1})^{2g}
+{\sqrt{-1}\over 4}(1-\sqrt{-1})^{2g}
=\sum_{j=0}^{g-1}\delta(2-j){2g\choose j}.$$
\end{pf}

\subsection{$\Theta$-characteristics}
Recall that the {\em $\Theta$-characteristics} are the fixed points of
the involution $\iota$ acting on $\Pic^{g-1}(X)$, i.e. they are the
(isomorphism classes of) line bundles $L$ on $X$ such that there is an 
isomorphism $L^{\otimes 2}\isomo \Omega^1_X$. A $\Theta$-characteristic
is called {\em even} (respectively {\em odd}) if $\dim H^0(X;L)$ is
even (respectively odd). The total number of $\Theta$-characteristics
on a curve of genus $g$ is equal to $2^{2g}$.

The fixed points of the action of $\iota$ on $\Theta$ correspond to the 
$\Theta$-characteristics $L$ with $\dim H^0(X;L)\geq 1$. Note that the fixed
point set $\Theta^\iota$ contains all of the odd $\Theta$-characteristics.

From Proposition \ref{prop:STr} we obtain the formula to the number of odd
$\Theta$-characteristics.
Naturally, this is a classical result, first proved
by Wirtinger in \cite{W} using theta functions.
An algebro-geometric proof valid in all odd characteristics
was given by Mumford in \cite{Mu1}.

\begin{prop}
The number of odd $\Theta$-characteristics is equal to
$2^{2g-1}-2^{g-1}$.
\end{prop}
\begin{pf}
The Lefschetz Fixed Point Formula (\cite{GM2}, \cite{V}) applied to $\iota$ 
gives
\[
STr(\iota, IH^\bullet(\Theta:\bbQ)) =
\sum_{L\in\Theta^\iota}STr(\iota, H^\bullet((i_{!*}\bbQ_{\Theta^{reg}})_L))
\]
where $(i_{!*}\bbQ_{\Theta^{reg}})_L)$ is the stalk of the the sheaf
$i_{!*}\bbQ_{\Theta^{reg}}$ at the point $L$.

By Proposition \ref{prop:IC-is-dir-im} there is an isomorphism
$H^\bullet((i_{!*}\bbQ_{\Theta^{reg}})_L)\isomo H^\bullet(\phi^{-1}(L);\bbQ)$,
and there is a natural identification $\phi^{-1}(L)\isomo\bbP(H^0(X;L))$.
Let $r(L) = \dim H^0(X;L) - 1$. Then we have
\[
H^p((i_{!*}\bbQ_{\Theta^{reg}})_L)\isomo\left\lbrace\begin{array}{lll}
\bbQ & \text{if} & 0\leq p=2j\leq r(L)\\ 0 & \text{otherwise}
\end{array}\right.
\]

By Lemma \ref{lemma:str-fp} below
\[
STr(\iota, H^\bullet((i_{!*}\bbQ_{\Theta^{reg}})_L)) =
\sum_{j=0}^{r(L)} (-1)^j =
\left\lbrace\begin{array}{lll} 1 & \text{if} & \text{$L$ is odd} \\
0 & \text{if} & \text{$L$ is even}
\end{array}\right.
\]
and the proposition follows immediately from Proposition \ref{prop:STr}.
\end{pf}

\begin{lemma}\label{lemma:str-fp}
\[
Tr(\iota : H^{2j}((i_{!*}\bbQ_{\Theta^{reg}})_L)\to
H^{2j}((i_{!*}\bbQ_{\Theta^{reg}})_L)) = (-1)^j\ .
\]
\end{lemma}
\begin{pf}
Consider the commutative diagram
\[
\begin{CD}
H^{2j}(S^{g-1}(X);\bbQ) @>>> H^{2j}(\phi^{-1}(L);\bbQ) \\
@V{\isomo}VV @V{\isomo}VV \\
IH^{2j}(\Theta;\bbQ) @>>> H^{2j}((i_{!*}\bbQ_{\Theta^{reg}})_L)
\end{CD}
\]
The top horizontal map is nonzero for $j=1$ since the class of an ample
divisor on $S^{g-1}(X)$ must have a nontrivial restriction. Therefore it
is nontrivial for all $j\leq r(L)$, hence surjective because the target
is one-dimensional. Since the vertical maps are isomorphisms so is
the bottom horizontal map.

Since the composition
\[
H^k(\Pic^{g-1}(X);\bbQ) @>>> H^k(S^{g-1}(X);\bbQ) @>>> H^k(\phi^{-1}(L);\bbQ)
\]
is trivial for all $k$ it follows that the compostion
\[
H^1(\Pic^{g-1}(X);\bbQ)\otimes IH^{2j-1}(\Theta) @>>>
IH^{2j}(\Theta;\bbQ) @>>> H^{2j}((i_{!*}\bbQ_{\Theta^{reg}})_L)
\]
must be trivial. Note that $\Coker(H^1(\Pic^{g-1}(X);\bbQ)
\otimes IH^{2j-1}(\Theta)\to IH^{2j}(\Theta;\bbQ))$ is generated by the image
of $1\otimes\eta_{g-1}^j$.

The induced map
\[
\Coker(H^1(\Pic^{g-1}(X);\bbQ)\otimes IH^{2j-1}(\Theta) @>>>
IH^{2j}(\Theta;\bbQ)) @>>> H^{2j}((i_{!*}\bbQ_{\Theta^{reg}})_L)
\]
is an isomorphism which is clearly $\iota$-equivariant and the lemma follows
from \eqref{formula:RR-eta}.
\end{pf}

\subsection{Un calcul encore plus triste}
Now we study the relation between the involution
$\iota$ acting on $IH^\bullet(\Theta;\bbQ)$ and the algebra
structure on $IH^\bullet(\Theta;\bbQ)$ induced by the isomorphism
$IH^\bullet(\Theta;\bbQ)\tilde{\to}H^\bullet(S^{g-1}(X);\bbQ)$.

\begin{prop}
For $g$ even, $\iota$ does not preserve
the algebra structure.
\end{prop}
\begin{pf}
According to \eqref{formula:RR-eta} we have $\iota(\eta_{g-1}) =
-\eta_{g-1} + \phi_{g-1}^*(\theta)$. Therefore the class $2\eta_{g-1}
+\phi_{g-1}^*(\theta)$ satisfies
\[
\iota(\eta_{g-1}+\phi_{g-1}^*(\theta)) = -(\eta_{g-1}
+\phi_{g-1}^*(\theta))\ .
\]

On the other hand, the class $(2\eta_{g-1} +\phi_{g-1}^*(\theta))^{g-1}$
satisfies
\[
\iota((\eta_{g-1}+\phi_{g-1}^*(\theta))^{g-1}) = (\eta_{g-1}
+\phi_{g-1}^*(\theta))^{g-1} = (-1)^{g-1}(\iota (\eta_{g-1}
+\phi_{g-1}^*(\theta))^{g-1})
\]
since $\iota$ acts trivially on the one dimensional space
$IH^{2g-2}(\Theta;\bbQ)$. It remains to observe that
$(\eta_{g-1} +\phi_{g-1}^*(\theta))^{g-1}\neq 0$.
\end{pf}

\end{document}